\documentclass[12pt]{article}
\usepackage{amsmath,amsfonts,amsthm,epsf}
\usepackage[colorlinks=true, pdfstartview=FitV, linkcolor=blue,
            menucolor=green, citecolor=blue, urlcolor=red]{hyperref}
\usepackage{graphicx,subfigure}
\usepackage{natbib}

\newtheorem{lemma}{Lemma}
\newtheorem{theorem}{Theorem}
\newtheorem{remark}{Remark}

\newcommand{\n}{n\,\nu_n}

\begin{document}

\title{Are adaptive allocation designs beneficial for improving power in binary response trials?
\author{
David Azriel$^1$, Micha Mandel$^1$ and Yosef Rinott$^{1,2}$\\ \\
$^1$ Department of statistics, The Hebrew University of Jerusalem, Israel \\ $^2$ Center for the Study of Rationality,
The Hebrew University of Jerusalem, Israel, \\and Luiss, Rome. \\
}
\thanks{We thank Amir Dembo for helpful comments and in particular for deriving the exact expression for $\nu^*$ in \eqref{eq:nu*}.} }
\maketitle

\begin{abstract}

We consider the classical problem of
selecting the best of two treatments in clinical trials with binary response. The target is to find the design that maximizes the power of the relevant test.
Many papers use a normal approximation to the power function and claim that \textit{Neyman allocation} that assigns subjects to treatment groups according to the ratio of the responses' standard deviations, should be used. As the standard deviations are unknown,  an adaptive
design is often recommended. The asymptotic justification of this approach is arguable, since it uses the normal approximation in tails where
the error in the approximation is larger than the estimated quantity. We consider two different approaches for optimality of designs that are related to Pitman and Bahadur definitions of relative efficiency of tests. We prove that the optimal allocation according to the Pitman criterion is the balanced allocation and that the optimal allocation according to the Bahadur approach depends on the unknown parameters.  Exact calculations reveal that the optimal allocation according to Bahadur is often close to the balanced design, and the powers of both are comparable to the Neyman allocation for small sample sizes and are generally better for large experiments. Our findings have important implications to the design of experiments, as the balanced design is proved to be optimal or close to optimal and the need for the complications involved in following an adaptive design for the purpose of increasing the power of tests is therefore questionable.
\end{abstract}

\vspace{0.5cm} KEYWORDS: Neyman allocation, adaptive design, asymptotic power, Normal approximation, Pitman efficiency, Bahadur efficiency, large deviations.

\section{Introduction}

We consider the problem of optimal allocation of individuals to two treatment groups with the goal of selecting the better treatment. The problem arises frequently in clinical trials, which  usually have several possibly conflicting purposes such as minimizing the number of subjects treated in the inferior treatment or maximizing the power of the relevant test. The current paper focuses on the latter goal and aims at answering the first question appearing in the chapter ``Fundamental questions of response-adaptive randomization'' of the book by \citet{Hu_Ros_2006}: what allocation maximizes power? It appears that the accepted answer to that question is the Neyman allocation, see references below. However, it is shown, both theoretically and by exact calculations, that the balanced allocation, that is, assigning an equal number of subjects to each treatment, is optimal or close to optimal. Unlike Neyman allocation, the balanced allocation does not depend on unknown parameters, and therefore no adaptive estimation is required. Adaptive designs are complex by nature, and our results question the need for the conducting such designs when the goal is to maximize power.

Let $A$ and $B$ be two treatments with unknown probabilities of
success, $p_A$ and $p_B$. A trial with $n$ subjects is planned with $N_A(n)$ and $N_B(n)$ subjects assigned to treatment $A$ and $B$, respectively, where $N_A(n)+N_B(n)=n$. For each subject, a binary response, success or failure, is observed.  Let $\nu_n := {N_A(n)}/n$ be the proportion of subjects assigned to treatment $A$. We sometimes refer to $\nu_n$ as the \textit{allocation}. The design problem considered here is of choosing the optimal allocation $\nu_n$ that maximizes the power of the standard test of the hypothesis $p_A = p_B$ versus one or two-sided alternatives. For given $n$,  $p_A$, and $p_B$, the optimal allocation fraction $\nu_n$ can be found by a finite search over all possible allocations. Here we study this problem for large $n$ instead, and look for the asymptotically optimal allocation fraction  $\nu^*$.


Let  $Y_i(m) \sim Bin(m,p_i)$ be the number of successes if $m$ patients are assigned to  treatment $i$ ($i=A,B$). Let also
${\hat p}_A={\hat p}_A(N_A(n))=\frac{Y_A(N_A(n))}{N_A(n)}$ and ${\hat p}_B={\hat p}_B(N_B(n))=\frac{Y_B(N_B(n))}{N_B(n)}$ be the estimators of $p_A$ and $p_B$;
note that ${\hat p}_A$ and ${\hat p}_B$ depend on $n$ and the allocation sequence $\nu_n$, however they are suppressed
for notational convenience.
The Neyman allocation rule, $\nu = \frac{\sqrt{p_A(1-p_A)}}{\sqrt{p_A(1-p_A)}+\sqrt{p_B(1-p_B)}}$, minimizes the variance of the estimator ${\hat p}_A(n) - {\hat p}_B(n)$ for the difference of probabilities (e.g., \cite{Melfi}). However, it is not clear that the Neyman allocation also maximizes the power of the Wald test for equality of proportions, as appears to be widely believed (e.g., \citet{Brittain}; \citet{Rosenberger};  \citet{Hu_Ros_2003}; \citet{Bandyopadhyay};  \citet{Hu_et_al}; \citet{Hu_Ros_2006};  \citet{Tymofyeyev}; \citet{Biswas}; \citet{Zhu}; \citet{Chambaz}). For example, when comparing the Neyman allocation to the balanced design, the latter authors claim that ``resorting to the balanced treatment mechanism may be a very poor (inefficient) choice''. Below we show that this claim is asymptotically incorrect.

The standard Wald statistic for comparing $p_A$ and $p_B$ is
\[
W:= \{\hat{p}_B - \hat{p}_A\}\Big/{\sqrt{V(\hat{p}_A,\hat{p}_B,n,\nu_n)}},
\]
where $V(p_A,p_B,n,\nu_n)=\frac{{p}_A(1-{p}_A)}{\nu_n \cdot n}+\frac{{p}_B(1-{ p}_B)}{(1-\nu_n)\cdot n}$.
In the above papers, the power is often calculated by approximating the distribution of the squared Wald statistic by a non-central chi-square distribution; the Neyman allocation then maximizes the non-centrality parameter.
The argument is based on the following normal approximation:
\begin{align*}
&P_{p_A,p_B}(W >z_{1-\alpha})=P_{p_A,p_B}\left(\frac{{\hat p}_B - {\hat p}_A - (p_B-p_A)}{\sqrt{V(p_A,p_B,n,\nu_n)}}>
\frac{z_{1-\alpha} \cdot {\sqrt{V(\hat{p}_A,\hat{p}_B,n,\nu_n)}} - (p_B-p_A)}{\sqrt{V(p_A,p_B,n,\nu_n)}}\right)\\
& \approx 1 - \Phi\left(\frac{z_{1-\alpha} \cdot {\sqrt{V(\hat{p}_A,\hat{p}_B,n,\nu_n)}} - (p_B-p_A)}{\sqrt{V(p_A,p_B,n,\nu_n)}}\right) \approx 1 - \Phi\left(z_{1-\alpha} - \frac{p_B-p_A}{{\sqrt{V(p_A,p_B,n,\nu_n)}}}\right),
\end{align*}
where $\Phi$ is the standard normal distribution function, and $z_{1-\alpha}=\Phi^{-1}(1-\alpha)$. The Normal approximation is  valid only if $({p_B-p_A})/{{\sqrt{V(p_A,p_B,n,\nu_n)}}}=O(1)$, i.e., when $p_B-p_A\approx n^{-1/2}$. However, for fixed $p_B-p_A > 0$, the term $({p_B-p_A})/{{\sqrt{V(p_A,p_B,n,\nu_n)}}}$ is of order $\sqrt{n}$, and the expression $\Phi\left(z_{1-\alpha} - \frac{p_B-p_A}{{\sqrt{V(p_A,p_B,n,\nu_n)}}}\right)$ is of asymptotic order that is smaller than the precision of the normal approximation, and
therefore its use is problematic. Thus, the claim that Neyman allocation maximizes the power seems theoretically questionable.

For asymptotic power comparisons and evaluation of the relative asymptotic efficiency of certain tests, two different criteria are often used, related to the notions of Pitman and Bahadur efficiency (see e.g., \citet{van der Vaart}, Chapter 14). In our context, the Pitman approach looks at sequences of probabilities $p_B^k > p_A^k$ that tend to a common limit at a suitable rate. The Bahadur approach considers fixed probabilities $p_A$ and $p_B$ and approximates the power using large deviations theory.

We show in the next sections that the optimal allocation corresponding to the Pitman approach is always $\nu^*=0.5$ while the Bahadur optimal allocation depends on $p_A$ and $p_B$ and can be calculated in a way described below. Interestingly, computation of the Bahadur criterion for different values of $p_A$ and $p_B$ reveals that the optimal allocation  is often close to $0.5$. In disagreement with some of the papers mentioned above, these results cast doubts on the asymptotic justification of  adaptive  designs and show that, at best, such designs can  lead to a  practically negligible improvement over a non-sequential balanced design in terms of power.

The paper is organized as follows: Sections \ref{sec:pitman} and \ref{sec:bahadur} describe the approaches of Pitman and Bahadur for maximizing the power, and find the corresponding optimal rules. In Section \ref{sec:calculation}, the optimal allocation according to the Bahadur criterion is calculated for different parameters and compared to the Neyman allocation. Exact calculations are performed for a wide range of parameters.  A related problem that arises in dose findings experiments is discussed in Section \ref{sec:diff}; the Neyman allocation is shown to be optimal or close to optimal in this case. Section \ref{sec:general} extends the Bahadur approach to general (rather than binary) responses;  concluding remarks are given in Section \ref{sec:conclusion}. All proofs are given in the Appendix.

\section{The Pitman Approach} \label{sec:pitman}

Pitman relative efficiency provides an asymptotic comparison of two families of tests applied to a sequence of statistical problems. Here we utilize the same idea to compare different allocation fractions.

Consider a sequence of statistical problems indexed by $k$, where $p^k_A = p + \frac{\delta_A }{\sqrt{k}}$, $p^k_B = p + \frac{\delta_B}{ \sqrt{k}}$,  for $\delta_A < \delta_B$ and $0<p<1$.
Let $n_k=n_k(\delta_A, \delta_B,p,\alpha,\beta,\{\nu_n\})$ be the minimal number of observations required for a one-sided Wald test at significance level $\alpha$ and power at least $\beta$ (for $\beta > \alpha$) at the point $p^k_A,p^k_B$, where the observations are allocated to the two groups according to the fraction $\nu_n$. Set $n_k=\infty$ if no finite number of observations satisfies these requirements. The next theorem implies that the balanced allocation is asymptotically optimal.

\begin{theorem}{\label{Pitman}}
Fix $\delta_A < \delta_B$, $\alpha<\beta$ and $0<p<1$. Let $\{\nu_n\}$ be a any sequence of allocations and let $\{\tilde{\nu}_n\}$ be another sequence of allocations satisfying $\tilde{\nu}_n\rightarrow 1/2$. Then
\[
 \liminf_{k \rightarrow \infty} \frac{n_k(\delta_A, \delta_B,p,\alpha,\beta,\{\nu_n\})}{n_k(\delta_A, \delta_B,p,\alpha,\beta,\{\tilde{\nu}_n\})} \ge 1.
 \]
\end{theorem}

The theorem follows readily from the following lemma, proved in the Appendix.

\begin{lemma}{\label{lem:Pitman}}
\begin{enumerate}
\renewcommand{\theenumi}{\Roman{enumi}}
\item If $\nu_n \rightarrow \nu$ for $0<\nu<1$ then
\begin{equation}\label{eq:Pitman}
\lim_{k \rightarrow \infty} \frac{n_k}{k} = \left(\frac{z_{1-\alpha} - z_{1-\beta}}{(\delta_B-\delta_A)}\cdot  {\sqrt{\frac{p(1-p)}{\nu(1-\nu)}}}\right)^2.
\end{equation}
\item If $\nu_n \rightarrow 0$ or $\nu_n \rightarrow 1$ then
\[
\lim_{k \rightarrow \infty} \frac{n_k}{k} = \infty
\]
\item For any sequence of allocations $\{\nu_n\}$
\begin{equation*}
\liminf_{k \rightarrow \infty} \frac{n_k}{k} \ge \left(\frac{z_{1-\alpha} - z_{1-\beta}}{(\delta_B-\delta_A)}\cdot  {\sqrt{\frac{p(1-p)}{\frac{1}{4}}}}\right)^2 .
\end{equation*}
\end{enumerate}
\end{lemma}

Theorem \ref{Pitman} holds also when considering a two-sided test.
The theorem shows that the balanced design is asymptotically optimal in the Pitman sense, and as a consequence, one cannot gain efficiency (in the above sense) by considering sequential adaptive designs.
The key point here is that when $p^k_A$ and $p^k_B$ converge to the same value $p$, the variances of their estimators converge to the same value and hence the limiting Neyman allocation is 1/2 regardless of $p$. This phenomenon is not observed in problems concerning the Normal distribution or similar cases where the variance is not a function of the mean.

It can be  argued that rather than considering sequences of statistical problems as above, one should optimize for fixed $p_A$ and $p_B$. The next section deals with this case.

\section {The Bahadur Approach} \label{sec:bahadur}


In this section, large deviations theory is used to approximate the power of the Wald test for fixed $p_A$ and $p_B$. This power increases exponentially to one with $n$ at a rate that depends on the allocation fraction $\nu$. Recall that
${\hat p}_A$ and ${\hat p}_B$ depend on both $n$ and an allocation $\nu_n$. The aim is to find the optimal limiting allocation fraction $\nu^*$ for which the rate is maximized.
We prove the following large deviations result:

\begin{theorem} \label{thm:test} Define
\[
H(t,\nu):= \nu \log(1-p_A+p_A e^{t/\nu})+(1-\nu)\log(1-p_B +p_B e^{-t/(1-\nu)}),
\]
and let ${g(\nu)}:= \inf_{t>0} H(t,\nu)$.
\begin{enumerate}
\renewcommand{\theenumi}{\Roman{enumi}}
\item One sided test: assume that $p_B > p_A$  and $\nu_n \rightarrow \nu$, where $0<\nu<1$, then for any constant $K \ge 0$
\begin{equation}\label{eq:one_sided}
\lim_n \frac{1}{n} \log\left\{1- P\left(\frac{{\hat p}_B -  {\hat p}_A}{\sqrt{V({\hat p}_A,{\hat p}_B,n,\nu_n)}}>K\right)\right\} = g(\nu).
\end{equation}
\item Two sided test: assume that $p_B \ne p_A$  and $\nu_n \rightarrow \nu$, where $0<\nu<1$, then for any constant $K>0$
\begin{equation}
\label{eq:two_sided1}
\lim_n \frac{1}{n} \log\left\{1- P\left(\frac{\{{\hat p}_B -  {\hat p}_A\}^2}{{V({\hat p}_A,{\hat p}_B,n,\nu_n)}}>K\right)\right\} = g(\nu).
\end{equation}
\item If $\nu_n \rightarrow 0$ or 1 then \eqref{eq:one_sided} and \eqref{eq:two_sided1} hold with $g(0)=g(1)=0$.
\end{enumerate}
\end{theorem}

Note that ${\hat p}_B - {\hat p}_A$ is not an average of $n$ i.i.d random variables and, therefore, Theorem \ref{thm:test} does not follow directly from the Cram\'{e}r-Chernoff theorem (see e.g., \citet{van der Vaart}, p. 205), however, its proof uses similar ideas.

For each fixed $n$, let $\nu^{*(1)}_n=\nu^{*(1)}_n(p_A,p_B,K)$ be the allocation that maximizes the power of the one sided test for a total sample size of $n$ subjects, i.e,
\[
\nu^{*(1)}_n={\arg\max}_{\nu_n \in \{\frac{1}{n},\ldots,\frac{n-1}{n}\}} P\left(\frac{{\hat p}_B -  {\hat p}_A}{\sqrt{V({\hat p}_A,{\hat p}_B,n,\nu_n)}}>K\right);
\]
similarly, $\nu^{*(2)}_n=\nu^{*(2)}_n(p_A,p_B,K)$ is the optimal allocation of the two-sided test.

Let $\nu^*=\nu^*(p_A,p_B) := {\arg \min}_\nu g(\nu)$. It is easy to prove directly that $g$ is strictly convex, and the minimum
is attained uniquely. More generally, it is readily shown by differentiation that if $M(t)=Ee^{tX}$ is a moment generating function, then $\nu M(t/\nu)$ is a convex function of $\nu$. Theorem \ref{thm:test} suggests the use of $\nu^*$ as the design fraction. However, for a given $n$, the optimal allocation, is not necessarily $\nu^*$, but the fraction $\nu^{*(1)}_n$ or $\nu^{*(2)}_n$ for the one or two-sided test, respectively. Therefore, it is reasonable to use $\nu^*$ as the design fraction only if $\nu^{* (i)}_n\rightarrow \nu^*$ for $i=1,2$. The following theorem shows that this is indeed the case.
\begin{theorem} \label{nu_star}
\begin{enumerate}
\renewcommand{\theenumi}{\Roman{enumi}}
\item If  $p_B > p_A$ then for any $K \ge 0$, $\nu^{*(1)}_n \rightarrow \nu^*$.
\item If  $p_B \ne p_A$ then for any $K > 0$, $\nu^{*(2)}_n \rightarrow \nu^*$.
\end{enumerate}
\end{theorem}

\begin{remark}
Another formulation of these results, for the one-sided case, say, is the following:
assume that $p_B > p_A$ then for any sequence $\nu_n$  and constant $K \ge 0$
\begin{equation*}
\liminf_n \frac{1}{n} \log\left\{1- P\left(\frac{{\hat p}_B -  {\hat p}_A}{\sqrt{V({\hat p}_A,{\hat p}_B,n,\nu_n)}}>K\right)\right\} \ge g(\nu^*),
\end{equation*}
and the infimum is attained for sequences $\nu_n \rightarrow \nu^*$.
\end{remark}

\begin{remark}
When ${p}_A<{p}_B$ represent success probabilities of two treatments, and treatment $B$ is selected as better if  ${\hat p}_B(n)> {\hat p}_A(n)$,
then  the expression in \eqref{eq:one_sided} with $K=0$  approximates the probability of incorrect selection.
\end{remark}

\section{Numerical Illustration} \label{sec:calculation}

Some tedious calculations show that
\begin{equation}\label{eq:nu*}
\nu^*=\log\left\{\frac{p_B\log(\frac{p_B}{p_A})}{(1-p_B)\log(\frac{1-p_A}{1-p_B})}\right\}\Big/\log\left\{\frac{p_B (1-p_A)}{p_A(1-p_B)}\right\} .
\end{equation}
Table \ref{tab1} compares the asymptotic Bahadur optimal allocation and the Neyman allocation for several pairs $(p_A, p_B)$. The table and further systematic numerical calculations indicate that the Bahadur allocation is closer to 0.5 than the Neyman allocation and that it is quite close to 0.5 unless $p_A$ and $p_B$ are very far apart (e.g., $p_A=0.5,p_B=0.9$). In the latter case, the power is close to 1 for any reasonable allocation. These findings justify the use of the balanced allocation and question the utility of more complicated adaptive sequential designs.

\begin{table} 
\caption{ \footnotesize The optimal Bahadur allocation $\nu^*$ for different parameters compared to Neyman allocation.}
  \begin{center}
  \begin{tabular}{cccc} \label{tab1}
    $p_A$ & $p_B$ & $\nu^*$ & Neyman allocation \\
     \hline \hline
    0.5 & 0.8 & 0.518 & 0.556 \\
    0.5 & 0.65 &  0.504 & 0.512 \\
    0.6 & 0.75 & 0.510 & 0.531 \\
    0.7 & 0.75 &  0.505 & 0.514 \\
    0.7 & 0.85 &  0.521 & 0.562 \\
    0.7 & 0.9  &  0.535  & 0.604 \\
    0.85 & 0.95 & 0.541 & 0.621 \\
    0.5 & 0.9 &   0.542 & 0.625 \\
  \end{tabular}
  \end{center}
\end{table}

\begin{figure}[ht!]
        \subfigure[Moderate sample size ($n=200$)]{%
            \includegraphics[width=0.5\textwidth]{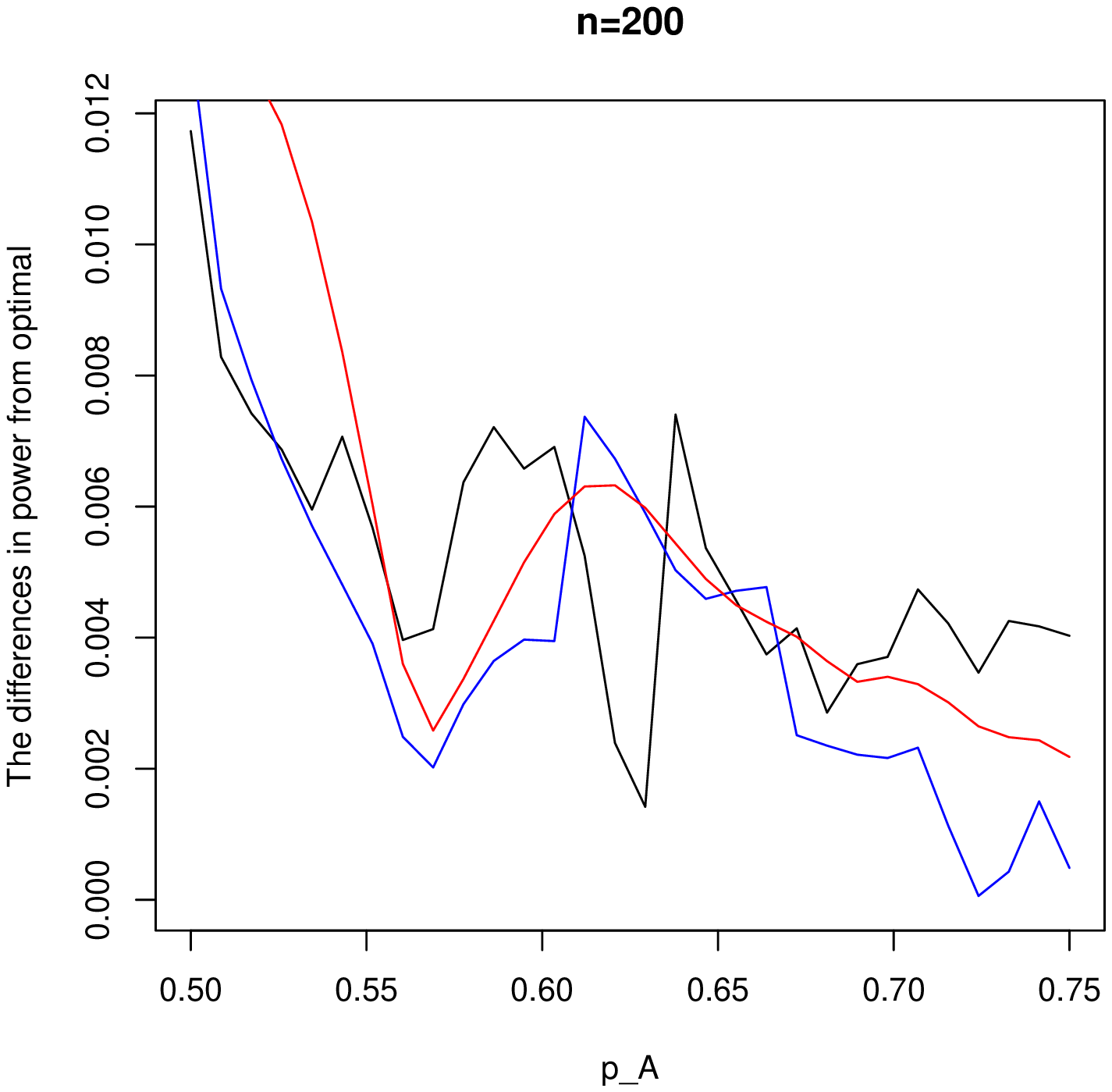}
        }%
        \subfigure[Large sample size ($n=500$)]{%
           \includegraphics[width=0.5\textwidth]{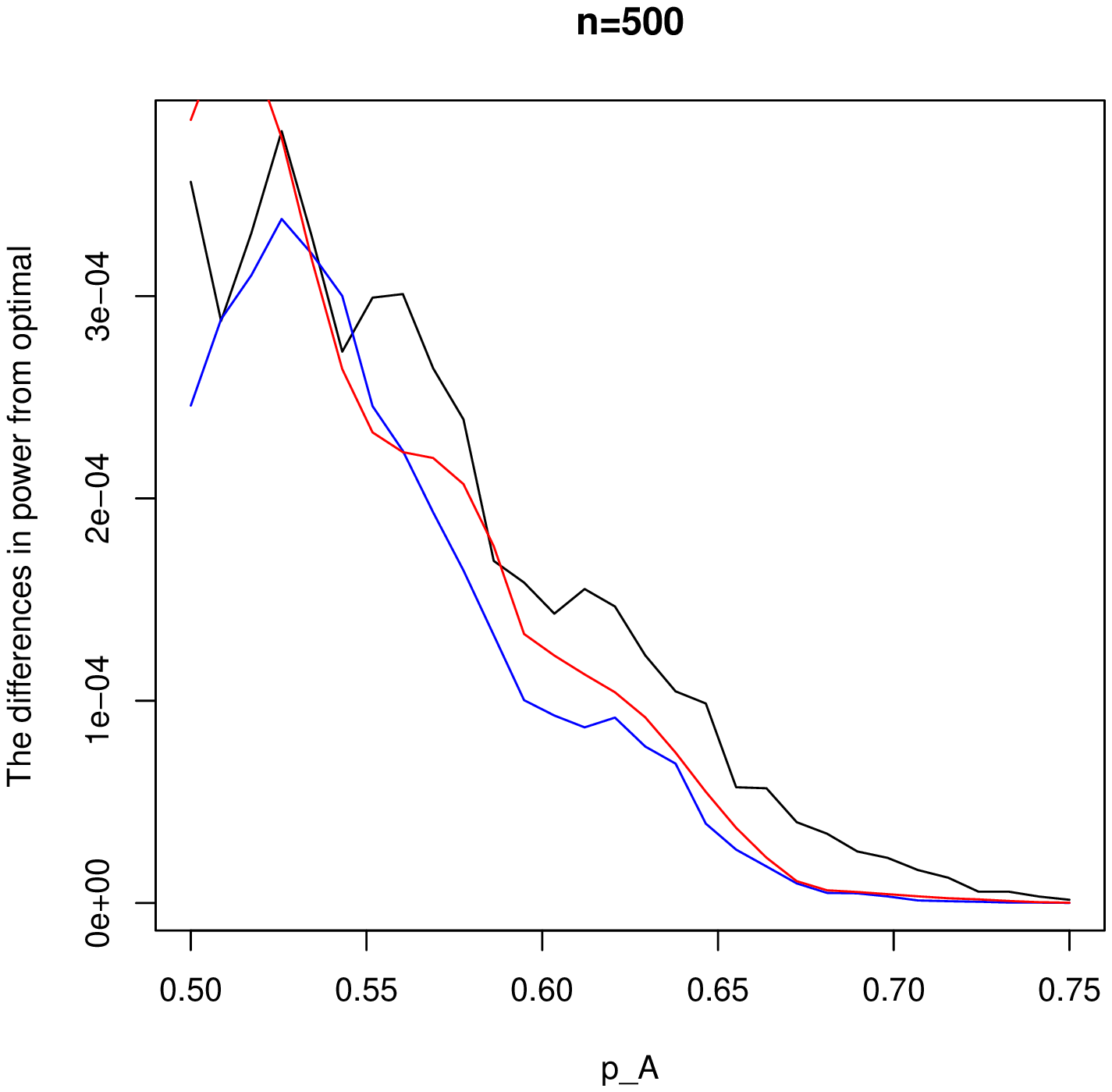}
        }\\ 
    \caption{\footnotesize{The differences between the maximal power of the two-sided test with critical value $K=1.96$, attainable (by $\nu^{*(2)}_n$) and the Neyman allocation (black), the balanced allocation (red) and the Bahadur allocation (blue) for $p_A=0.5,\ldots,0.75$, $p_B = p_A+0.2$; for moderate ($n=200$) and large ($n=500$) sample size.}}\label{fig:fig1}
\end{figure}

We preformed some exact calculations to compare the Bahadur allocation, the balanced allocation and the Neyman allocation. Figure  \ref{fig:fig1} compares the difference between the maximal possible power for sample size 200 and 500, and the power under the different allocation methods for the two-sided test with $\alpha=0.05$ and for different parameters. The power is calculated exactly using R. While for moderate sample size ($n=200$) no allocation is better for all the parameters we considered, for large sample size ($n=500$), Bahadur is better for almost all parameters, and the balanced allocation is usually better than Neyman; however, the differences in power are relatively small.

\begin{figure}[h!]
\centering
\includegraphics[height=9cm]{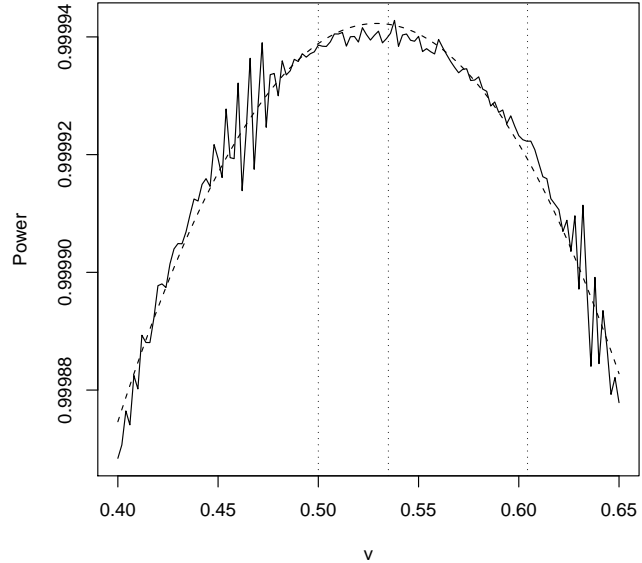}
\caption{\footnotesize The power of Wald tests with critical value $K=1.96$ for different allocations $\nu$ where $p_A=0.7 , p_B = 0.9$  and $n=500$. The smooth line is a parabolic fit to the function. The vertical lines show the balanced allocation ($\nu = 0.5$), the Bahadur allocation ($\nu =0.5349374 $) and the Neyman allocation ($\nu =  0.6043561$).  \label{fig:fig3}}
\end{figure}

Figure \ref{fig:fig3} shows the power of the two-sided test for different allocations where $p_A=0.7 , p_B = 0.9$; it is clearly seen that the Neyman allocation, which is widely recommended for maximizing the power, is far from being optimal. Thus, the exact calculations presented in this section support the theoretical results: the balanced allocation is usually better than the Neyman allocation for large samples, and they are indistinguishable for small samples. In all cases, the differences are quite negligible, and therefore the balanced allocation should be preferred due to its simplicity.

\section{A Related Problem}\label{sec:diff}

Dose finding studies are conducted as part of phase I clinical trials in order to find the maximal tolerated dose (MTD) among a finite, usually very small, number of potential doses. The MTD is defined as the dose with the closest  probability of toxic reaction to a pre-specified probability $p_0$.
Recently, we showed that under certain natural assumptions, in order to estimate the desired dose consistently, one
can consider experiments that eventually concentrate on two doses (\citet{Azriel}). Thus, asymptotically, the allocation problem in MTD studies reduces to the problem of finding which of two probabilities of toxic reaction $p_A<p_B$ (corresponding to the doses $d_A<d_B$) is closer to $p_0$.

Let $\hat{p}_A$ and ${\hat p}_B$ denote the proportions of toxic reactions in doses $d_A$ and $d_B$ based on a total sample size of $n$ individuals and an allocation $\nu_n$. For large $n$,  $\hat{p}_A<{\hat p}_B$, and a natural estimator for the MTD is $\widehat{MTD}=d_A$ if
$(\hat{p}_A+{\hat p}_B)/2 > p_0$ and $\widehat{MTD}=d_B$ otherwise. Similar to the problems discussed in previous sections, an
optimal design is an allocation rule of $n \, \nu_n$ and $n (1- \nu_n)$ individuals
to doses $d_A$  and $d_B$, respectively, such that
$P(\widehat{MTD}=d_A)=P((\hat{p}_A+{\hat p}_B)/2 > p_0)$ is maximized if $d_A$ is indeed the MTD.

For the current problem, the Pitman approach is translated to a comparison of designs under sequences of parameters $p^k_A$, $p^k_B$ and $p_0^{k}$ such that
$|(p^k_A+p^k_B)/2 - p_0^{k}|  = K/\sqrt{k}$, for fixed $0<K<\infty$, and $p^k_A \rightarrow p_A$,  $p^k_B \rightarrow p_B$.
Let $0<\nu<1$ and let $n_k=n_k(p^k_A, p^k_B, p_0^{k},\alpha,\{\nu_n\})$ be the minimal number of observations required such that the probability of incorrect estimation of the MTD is smaller than $\alpha$ for the given parameters when the allocation  for dose $d_A$ is $n\cdot\nu_n$. As in Lemma \ref{lem:Pitman}, it can be shown that if $\nu_n \rightarrow \nu$ then
\[
\lim_{k \rightarrow \infty} \frac{n_k}{k} = \left\{ \frac{z_{1-\alpha}}{2K} \right\} ^2 \left\{\frac{p_A(1-p_A)}{\nu} + \frac{p_B(1-p_B)}{1-\nu}\right\}.
\]
Thus,
the asymptotically optimal design uses Neyman allocation,  $\nu = \frac{\sqrt {p_A(1-p_A)}}{\sqrt {p_A(1-p_A)}+\sqrt {p_B(1-p_B)}}$, as it minimizes the limit of ${n_k}/{k}$. Unlike the previous problem, now $p^k_A$ and $p^k_B$ do not converge to the same value under the Pitman approach as defined here, and hence the Neyman allocation does not reduce to the balanced design.

For the case of fixed $p_A$, $p_B$, and $p_0$, assume that $p_B$ is  nearer than $p_A$ to
$p_0$, and consider the problem of minimizing the probability of selecting $d_A$. The following theorem, analogous to Theorems \ref{thm:test} and \ref{nu_star}, gives the asymptotic optimal allocation rule in the current setting.
\begin{theorem} \label{large_deviation}
Let $\nu_n=N_A(n)/n$, $0<\nu<1$,  and assume that
 $\nu_n \rightarrow \nu $, then,
\[
\lim_{n \rightarrow \infty} \frac{1}{n} \log P[\{{\hat p}_A+{\hat p}_B\}/2 \geq p_0] = \psi(\nu),
\]
where $\psi(\nu)=\inf_t\{\nu \log(1-p_A+p_A
e^{t/\nu})+(1-\nu)\log(1-p_B+p_B e^ {t/(1-\nu)})- 2 p_0t\}$.

Moreover, let $\nu^*=\arg \min \psi(\nu)$, and let $\nu^*_n$ be the value of
the allocation minimizing $P[\{{\hat p}_A+{\hat p}_B\}/2 \geq p_0]$  for a given $n$. Then, $\nu^*_n\rightarrow \nu^*$.
\end{theorem}

\begin{table} 
\caption{ \footnotesize Comparison of the Bahadur and Pitman (here Neyman) allocation rules for different parameters.}
  \begin{center}
  \begin{tabular}{ccccc} \label{tab2}
    $p_A$ & $p_B$ & $p_0$ & Bahadur & Pitman \\
     \hline \hline
    0.1 & 0.3 &0.28 & 0.420  & 0.396 \\
    0.2 & 0.35 & 0.3 & 0.460  & 0.456 \\
    0.22 & 0.33 & 0.3 & 0.471 & 0.468 \\
    0.25 & 0.35 & 0.33 & 0.479 & 0.476 \\
    0.2 & 0.4 & 0.33 & 0.455  & 0.449\\
    0.1 & 0.4 & 0.3  &  0.400  & 0.380 \\
  \end{tabular}
  \end{center}
\end{table}

We calculated $\nu^*$ for several values of $p_A$ and $p_B$ and found that $\nu^*$ is often close to the Neyman allocation, see Table \ref{tab2}.
Both criteria, Bahadur and Pitman, yield quite similar results in this problem. Allocating subjects according to the Neyman or Bahadur improves the probability of correct MTD estimation compared to the balanced allocation for very large samples, as the optimal allocations according to Bahadur or Pitman are far from 0.5. Calculations not presented here, show that for practical sample sizes for the MTD problem, all three methods differ in a negligible way.

\section{A General Response}\label{sec:general}

In previous sections, we dealt with the very important, though specific, case of a binary response. In this section, we consider the more general case where the response of an individual treated in group A (B) follows a distribution $F_A$ ($F_B$) having moment generating function $M_A(t)$ ($M_B(t)$), and find the optimal allocation according to the Bahadur approach. Let $\bar Y_A(m)$ ($\bar Y_B(m)$) denote the average of $m$ responses of subjects having treatment $A$ ($B$). Assume that the treatment with the largest mean response is declared better at the end of the experiment. The following theorem, which can be proved in a similar way  as Theorems \ref{thm:test} and \ref{nu_star}, provides the Bahadur optimal allocation rule for correct selection:
\begin{theorem} \label{large_deviation_gen}
Assume that treatment $B$ is better, i.e, $\int x  F_B(dx) > \int x  F_A(dx)$,  and that
 $\nu_n \rightarrow \nu $. Then,
\[
\lim_{n \rightarrow \infty} \frac{1}{n} \log P\left\{\bar{Y}_A(n\,\nu_n) \ge \bar{Y}_B(n(1-\nu_n))\right \} = h(\nu),
\]
where
\begin{equation} \label{eq:optimal_general}
h(\nu)=\inf_t[\nu \log\{M_A(t/\nu)\}+(1-\nu)\log\{M_B[-t/(1-\nu)]\}].
\end{equation}
Moreover, let $\nu^*=\arg \min_\nu h(\nu)$, and $\nu^*_n$ be the value of
the allocation minimizing \\
$P\left\{\bar{Y}_A(n\,\nu_n) \ge \bar{Y}_B(n(1-\nu_n))\right \}$. Then $\nu^*_n\rightarrow \nu^*$.
\end{theorem}


When the responses in the two treatments are normally distributed, then the Bahadur allocation agrees with the Neyman allocation. This can be  easily verified by using the moment generating functions of Normal variables in (\ref{eq:optimal_general}). However, for other distributions, the allocations suggested by the Bahadur and the Neyman criteria may differ considerably. Table \ref{tab3} compares the Bahadur  and the Neyman allocations for different Poisson and Gamma distributions. The two rules clearly differ. As in the Binomial case, the Bahadur allocation is closer to 0.5 than to the Neyman allocation. Further study is required to determine if the improvement over the balanced allocation, in terms of power or probability of correct selection, is significant. Anyway, optimality of the Neyman allocation  for non-normal distributions should be questioned, and may hold only under restrictive
conditions.

\begin{table} 
\caption{ \footnotesize The optimal allocation $\nu^*$ for different distributions compared to the Neyman allocation.}
  \begin{center}
  \begin{tabular}{cccc} \label{tab3}
    $F_A$ & $F_B$  & Bahadur allocation & Neyman allocation \\
     \hline \hline
    Poisson(1) & Poisson(2) &  0.471 & 0.414 \\
    Poisson(2) & Poisson(3) & 0.483 & 0.449 \\
     Poisson(3) & Poisson(4) & 0.488 & 0.464 \\
    Poisson(4) & Poisson(5) &  0.491 & 0.472 \\
    Gamma(0.5,0.5) & Gamma(0.5,0.6) &0.515 &0.590 \\
    Gamma(0.5,0.5) & Gamma(0.5,0.7) &  0.528 &0.662 \\
    Gamma(0.5,0.5) & Gamma(0.5,0.8) & 0.539 & 0.719 \\
    Gamma(0.5,0.5) & Gamma(0.5,0.9) & 0.549 & 0.764 \\
  \end{tabular}
  \end{center}
\end{table}

\section{Conclusions}\label{sec:conclusion}

We discussed asymptotic approximations of power and probability of correct selection in testing and
selecting the best treatment, and in MTD finding, and related optimal allocation of subjects to treatments.

Neyman allocation is optimal when the response is Normal, and it is asymptotically optimal in the Pitman sense, that is, for converging sequences of alternatives as described above. In the binary response selection problem in which $p_A$ and
$p_B$ become closer, Neyman allocation reduces to a balanced allocation, independent of the parameters $p_A$ and
$p_B$. The Bahadur allocation for fixed $p_A$ and
$p_B$ turns out to be close to balanced, and therefore, by both criteria, our conclusion is that for the purpose of maximizing the power, adaptive allocation seems unwarranted, and the simpler, non-sequential balanced allocation should be preferred. However, when other criteria (e.g., ethical criteria) are of primary concern, as is often the case in clinical trials, the balanced design is not necessarily optimal and adaptive designs may be found beneficial.

Our findings are partly in contrast with the literature that bases allocations on the noncentrality parameter
appearing in a Normal or Chi-Square approximation (e.g., \citet{Rosenberger};  \citet{Tymofyeyev}).
These designs minimize or control the variance of the difference but need not be efficient in the sense of controlling or maximizing the power.


\appendix

\section{Appendix}\label{sec:proofs}

\noindent{\bf Proof of Lemma \ref{lem:Pitman} part I.}
The proof, included here for completeness,   uses arguments as in Theorem 14.19 in  \citet{van der Vaart} (p. 205), which is stated in terms of  relative efficiency rather than allocation.

First note  that $\lim_k n_k=\infty$; otherwise, there exists a bounded subsequence of $n_k$ on which the power  converges to
a value $\le \alpha$, since as $k \rightarrow \infty$ we have
$p^k_A-p^k_B  \rightarrow 0$. This contradicts the definition of $n_k$ and the assumption that $\alpha < \beta$.

By the Berry-Esseen theorem we have
\[
\frac{{\hat p}^k_A-(p+\frac{\delta_A}{\sqrt{k}})}{\sqrt{\frac{(p+\frac{\delta_A}{\sqrt{k}})[1-(p+\frac{\delta_A}{\sqrt{k}})]}{\nu_{n_k} \cdot n_k}}} \stackrel{D}{\rightarrow} N(0,1)
\]
since the third moment is bounded;  a similar limit holds for ${\hat p}^k_B$. Here we use the notation ${\hat p}_A^k={\hat p}_A(\nu_{n_k} \, n_k)=Y^k_A(\nu_{n_k}\, n_k)/(\nu_{n_k} \, n_k)$,
where $Y^k_A(m) \sim Bin(m, p^k_A)$ is the sum of $m$ independent binary responses with probability
$p^k_A$.

Now, if $\nu_{n_k} \rightarrow \nu$ we have
\begin{equation}\label{eq:CLT}
U_k:=\frac{\sqrt{n_k}( {\hat p}^k_B- {\hat p}^k_A)-(\delta_B-\delta_A)\sqrt{\frac{n_k}{k}}}{\sqrt{\frac{p(1-p)}{\nu(1-\nu)}}}
 \stackrel{D}{\rightarrow} N(0,1).
\end{equation}
Since $n_k \rightarrow \infty$, the critical value for the level $\alpha$ one-sided Wald test is $z_{1-\alpha}+o(1)$; then
\begin{align*}
&P_{p^k_A,p^k_B}(W >z_{1-\alpha}+o(1))= P\left(\frac{{\hat p}^k_B - {\hat p}^k_A}{\sqrt{V({\hat p}^k_A,{\hat p}^k_B,n_k,\nu_{n_k})}}>
z_{1-\alpha}+o(1)\right)\\
&=P\left(U_k >
\frac{(z_{1-\alpha}+o(1))\sqrt{V({\hat p}^k_A,{\hat p}^k_B,n_k,\nu_{n_k})n_k}}{\sqrt{\frac{p(1-p)}{\nu(1-\nu)}}}-\frac{(\delta_B-\delta_A)\sqrt{\frac{n_k}{k}}}{\sqrt{\frac{p(1-p)}{\nu(1-\nu)}}}\right).
\end{align*}
Also,
\[
\frac{\sqrt{V({\hat p}^k_A,{\hat p}^k_B,n_k,\nu_{n_k})n_k}}{\sqrt{\frac{p(1-p)}{\nu(1-\nu)}}}\stackrel{a.s.}{\rightarrow} 1,
\]
and since the limiting power is exactly $\beta$ we have due to (\ref{eq:CLT})
\[
z_{1-\alpha}-\frac{(\delta_B-\delta_A)\{\lim_k \sqrt{\frac{n_k}{k}}\}}{\sqrt{\frac{p(1-p)}{\nu(1-\nu)}}} = z_{1-\beta};
\]
hence \eqref{eq:Pitman} holds.

\noindent{\bf Proof of part II.}
We only prove the case  $\nu_n \rightarrow 0$, as $\nu_n \rightarrow 1$ is similar. If $n\nu_n$ is bounded, then the power converges to $\alpha$  and $n_k=\infty$ for large $k$.

Assume now that $n\,\nu_n \rightarrow \infty$; by the Berry-Esseen theorem and Slutsky's Lemma we have
\[
\sqrt{n_k\,\nu_{n_k}}({\hat p}^k_A-(p+\frac{\delta_A}{\sqrt{k}})) \stackrel{D}{\rightarrow} N(0,p(1-p))\,~\text{and}~\,\sqrt{n_k\,\nu_{n_k}}({\hat p}^k_B-(p+\frac{\delta_B}{\sqrt{k}})) \stackrel{D}{\rightarrow} 0.
\]
This implies that
\[
\sqrt{n_k\,\nu_{n_k}}({\hat p}^k_B-{\hat p}^k_A)-(\delta_B-\delta_A)\sqrt{\frac{n_k\,\nu_{n_k}}{k}} \stackrel{D}{\rightarrow} N(0,p(1-p))
\]
and by arguments as in the first part we have
\[
\lim_k \frac{n_k\,\nu_{n_k}}{k} = \left( \frac{z_{1-\alpha}-z_{1-\beta}}{\delta_B-\delta_A} \right)^2 p(1-p).
\]
Because $\nu_{n_k} \rightarrow 0$, $\lim_k \frac{n_k}{k} = \infty$.

\noindent{\bf Proof of part III.} There exists a subsequence $\{k'\}$ such that  $\nu_{n_{k'}} \rightarrow \nu'$ for some $\nu'$ and
\[
\liminf_{k \rightarrow \infty} \frac{n_{k}}{k} = \lim_{k' \rightarrow \infty} \frac{n_{k'}}{k'} = \left(\frac{z_{1-\alpha} - z_{1-\beta}}{(\delta_B-\delta_A)}\cdot  {\sqrt{\frac{p(1-p)}{\nu'(1-\nu')}}}\right)^2
,\]
where the second equality follows by part I. If $\nu'(1-\nu')=0$ we interpret the limit as $\infty$;  since $\nu'(1-\nu')\le \frac{1}{4}$ the third part of the lemma follows. \qed

\noindent {\bf Proof of Theorem \ref{thm:test} parts I and II.}
The proof follows known large deviations ideas; however, certain variations are needed for the present non-standard case. Notice that the probability in part I is larger than the probability of part II (for $\sqrt{K}$). Therefore, it is enough to show that for any $K \ge 0$
\begin{equation}\label{eq:one_sided1}
\limsup_n \frac{1}{n} \log\left\{1- P\left(\frac{{\hat p}_B -  {\hat p}_A}{\sqrt{V({\hat p}_A,{\hat p}_B,n,\nu_n)}}>K\right)\right\} \le g(\nu),
\end{equation}
and for any $K>0$
\begin{equation*}
\liminf_n \frac{1}{n} \log\left\{1- P\left(\frac{\{{\hat p}_B -  {\hat p}_A\}^2}{{V({\hat p}_A,{\hat p}_B,n,\nu_n)}}>K\right)\right\} \ge g(\nu).
\end{equation*}
In fact,  instead of the latter inequality we prove in the sequel a stronger result, namely
\begin{equation}
\label{eq:two_sided}
\liminf_n \frac{1}{n} \log P\left(0 \le \frac{{\hat p}_A -  {\hat p}_B}{\sqrt{V({\hat p}_A,{\hat p}_B,n,\nu_n)}}\le K'\right) \ge g(\nu),
\end{equation}
for all $K'>0$, which is also used for the case of $K=0$ in part I, when $K'=\infty$.

For the upper bound \eqref{eq:one_sided1},
define $S(n):=\sqrt{V({\hat p}_A,{\hat p}_B,\nu_n,n) \cdot n}$; notice that $S(n)$ is bounded.
Hence, for any $\varepsilon > 0$ and for large enough $n$ we have
\[
1- P\left(\frac{{\hat p}_B -  {\hat p}_A}{\sqrt{V({\hat p}_A,{\hat p}_B,n,\nu_n)}}>K\right) =
P\left({\hat p}_A -  {\hat p}_B \geq -\frac{K}{\sqrt n} S(n) \right) \leq P\left({\hat p}_A -  {\hat p}_B \geq - \varepsilon\right).
\]
Now, for any $t>0$,
\[
 P({\hat p}_A -  {\hat p}_B \geq - \varepsilon) = P(e^{t(\frac{Y_A(n\, \nu_n)}{N_A(n)/n} -  \frac{Y_B(n (1-\nu_n))}{N_B(n)/n})} \geq e^{-n t \varepsilon}) \leq E[e^{t(\frac{Y_A(n\,\nu_n)}{\nu_n} -  \frac{Y_B(n (1-\nu_n))}{1-\nu_n})}]e^{n t \varepsilon},
\]
by Markov's inequality. We can write the latter term as
\[
(1-p_A + p_A e^{t/\nu_n})^{\n}\cdot(1-p_B + p_B e^{-t/(1-\nu_n)})^{n (1-\nu_n)}e^{n t \varepsilon}
.\]
Since $\nu_n \rightarrow \nu$, and the inequality holds for all $t>0$,
\begin{equation*}
\limsup_n\frac{1}{n} \log P\left ({\hat p}_A -  {\hat p}_B \geq -\frac{K}{\sqrt n} S(n)\right) \le g_\varepsilon(\nu),
\end{equation*}
where $g_\varepsilon(\nu):= \inf_{t>0}\{\varepsilon t+ H(t,\nu)\}$.
This is true for any $\varepsilon >0$, and by the continuity of $g_\varepsilon(\nu)$ in $\varepsilon$ we have for any $K\ge 0$
\begin{equation*} 
\limsup_n\frac{1}{n} \log P\left({\hat p}_A -  {\hat p}_B \geq -\frac{K}{\sqrt n} S(n)\right) \le g(\nu),
\end{equation*}
which verifies \eqref{eq:one_sided1}.

To prove \eqref{eq:two_sided}, assume without loss of generality that $p_B > p_A$; define
\[
T_n:={\hat p}_A(n\,\nu_n) - {\hat p}_B(n (1-\nu_n) ) = \frac{Y_A(n\,\nu_n)}{n\,\nu_n}-\frac{Y_B(n(1-\nu_n))}{n(1-\nu_n)}.
\]
The log of the moment generating function of $T_n$ is
\begin{equation}\label{eq:H}
\log E[e^{t T_n}] = n\, \nu_n \log(1-p_A + p_A e^{\frac{t}{n\,\nu_n}})+n(1-\nu_n)\log(1-p_B + p_B e^{\frac{-t}{n(1-\nu_n)}}) = n H(\frac{t}{n},\nu_n).
\end{equation}
Since  $E[T_n]=p_A - p_B<0$,  by \eqref{eq:H}  we have $\frac{d}{dt}H(0,\nu_n)<0$. Also, $H(0,\nu_n)=0$ and $H(\cdot,\nu_n)$ is strictly convex being the log of a moment generating function, up to a constant. Since $P(T_n>0)>0$ it follows that $H(t,\nu_n)\rightarrow \infty$ as $t \rightarrow \infty$ and therefore, $\arg \min_{t>0} H(t,\nu_n) =: t_0^{(n)}$ is a unique interior point and $\frac{\partial}{\partial t}H(t_0^{(n)},\nu_n)=0$. Let $t_0$ be the minimizer of $H(\cdot,\nu)$; we show that $t_0^{(n)}\rightarrow t_0$. If there is a subsequence $\{t_0^{(n_k)}\}$ that converges to $t_1 \le \infty$ then $H(t_0^{(n_k)},\nu_{n_k}) \le H(t_0,\nu_{n_k})$ (as $t_0^{(n_k)}$ is the minimizer) implies $H(t_1,\nu) \le H(t_0,\nu)$ and therefore $t_1=t_0$ as the minimizer is unique and finite.

Define a new random variable $Z_n$, which is the Cram\'{e}r transform of $T_n$
\[
P(Z_n=z):= e^{-n g(\nu_n)}e^{z t_0^{(n)} n}P(T_n=z).
\]
Now,
\begin{align*}
&P\left(0 \le \frac{{\hat p}_A -  {\hat p}_B}{\sqrt{V({\hat p}_A,{\hat p}_B,n,\nu_n)}}\le K\right) = P\left(0\le T_n  \le \frac{K}{\sqrt{n}} S(n)\right)\\
&= E[ I\{0\le Z_n  \le \frac{K}{\sqrt{n}} S(n)\}e^{-Z_n t_0^{(n)} n}]e^{n g(\nu_n)} \ge P\big(0\le Z_n  \le \frac{K}{\sqrt{n}} S(n)\big) e^{-K \frac{1}{2}\sqrt{\frac{1}{\nu_n}+\frac{1}{1-\nu_n}} t_0^{(n)} \sqrt{n}}e^{n g(\nu_n)},
\end{align*}
where the last inequality holds since $e^{-Z_n} \ge e^{-\frac{K}{\sqrt{n}}S(n)} \ge e^{-\frac{K}{\sqrt{n}}\frac{1}{2}\sqrt{\frac{1}{\nu_n}+\frac{1}{1-\nu_n}}}$. It follows that
\[
g(\nu_n)-\frac{1}{n} \log P\left(0 \le \frac{{\hat p}_A -  {\hat p}_B}{\sqrt{V({\hat p}_A,{\hat p}_B,n,\nu_n)}}\le K\right) \le \frac{ -\log P\big(0\le Z_n  \le \frac{K}{\sqrt{n}} S(n)\big)}{n} + \frac{K \frac{1}{2}\sqrt{\frac{1}{\nu_n}+\frac{1}{1-\nu_n}} t_0^{(n)}}{\sqrt{n}}.
\]
Clearly, the second term on the right-hand side vanishes as $n$ goes to infinity; for the first, we claim that $\sqrt{n} Z_n$ is asymptotically $N(0,\frac{\partial^2}{\partial t^2}H(t^{(n)}_0,\nu_n))$ and consequently $P\big(0\le Z_n  \le \frac{K}{\sqrt{n}} S(n)\big)\rightarrow C$ for some constant $C>0$. Indeed, the log of the moment generating function of $\sqrt{n} Z_n$ is
\[
\log E[e^{s \sqrt{n} Z_n}]= -n g(\nu_n) + \log E[e^{T_n(s \sqrt{n}+ t_0^{(n)}n) }] = n\{ -H(t^{(n)}_0,\nu_n)+H(t_0^{(n)} + \frac{s}{\sqrt{n}},\nu_n)\},
\]
where the last equality follows from \eqref{eq:H} and the identity $g(\nu_n) =   H(t_0^{(n)},\nu_n)$. By Taylor expansion of $H(\cdot,\nu_n)$ around $t^{(n)}_0$ we obtain
\[
H(t_0^{(n)} + \frac{s}{\sqrt{n}},\nu_n)-H(t^{(n)}_0,\nu_n)=\frac{1}{2}\frac{s^2}{n}\frac{\partial^2}{\partial t^2}H(t^{(n)}_0,\nu_n) + O(n^{-3/2})
\]
since the first derivative is 0, and therefore,
\[
\log E[e^{s \sqrt{n} Z_n}] \rightarrow \frac{s^2}{2} \frac{\partial^2}{\partial t^2}H(t^{(n)}_0,\nu_n).
\]
We conclude that
\[
\limsup_n \left\{g(\nu_n)-\frac{1}{n} \log P\left(0 \le \frac{{\hat p}_A -  {\hat p}_B}{\sqrt{V({\hat p}_A,{\hat p}_B,n,\nu_n)}}\le K\right) \right\} \le 0,
\]
hence,
\[
\liminf \frac{1}{n} \log P\left(0 \le \frac{{\hat p}_A -  {\hat p}_B}{\sqrt{V({\hat p}_A,{\hat p}_B,n,\nu_n)}}\le K\right) \ge g(\nu)
\]
and part I and II follow.

\noindent {\bf Proof of part III.}
First note that \eqref{eq:one_sided1} clearly holds with $g(\nu)=0$ as $\log \{1-P(\cdot) \}\le 0$, so it remains
to prove \eqref{eq:two_sided} for $g(\nu)=0$,  that is, for any $K >0$
\begin{equation*}
\liminf_n \frac{1}{n} \log P\big( 0 \le  \frac{ {\hat p}_A -  {\hat p}_B}{\sqrt{V({\hat p}_A,{\hat p}_B,\nu_n,n)}} \leq K \big) \ge 0.
\end{equation*}
We only prove the case  $\nu_n \rightarrow 0$, as $\nu_n \rightarrow 1$ is similar.  If $ n\,\nu_n \not \rightarrow  \infty$  then ${\hat p}_A$ is inconsistent and the limit is easily seen to be zero. Assume now that $n\,\nu_n \rightarrow \infty$;
since
\[
V({\hat p}_A,{\hat p}_B,\nu_n,n) = \frac{{\hat p}_A (1-{\hat p}_A)}{n\,\nu_n} + \frac{{\hat p}_B (1-{\hat p}_B)}{n(1-\nu_n)} \ge  \frac{{\hat p}_A (1-{\hat p}_A)}{n\,\nu_n}
\]
we have
\[
 P\big( 0 \le  \frac{ {\hat p}_A -  {\hat p}_B}{\sqrt{V({\hat p}_A,{\hat p}_B,\nu_n,n)}} \leq K \big) \ge  P\big( 0 \le {\hat p}_A -  {\hat p}_B \le  \frac{K\sqrt{{\hat p}_A (1-{\hat p}_A)}}{\sqrt{n\, \nu_n}}  \big).
\]
Now, for $\varepsilon:= \frac{K p_A (1-p_A)}{2}$,
\begin{align*}
&P\big( 0 \le {\hat p}_A -  {\hat p}_B \le  \frac{K\sqrt{{\hat p}_A (1-{\hat p}_A)}}{\sqrt{n\, \nu_n}}  \big)\ge P\big(\{ 0 \le {\hat p}_A -  {\hat p}_B \le  \frac{K\sqrt{{\hat p}_A (1-{\hat p}_A)}}{\sqrt{n\, \nu_n}} \} \cap \{{\hat p}_B \in ( p_B - \frac{\varepsilon}{\sqrt{n\,\nu_n}}, p_B) \}    \big)  \\
&\ge P\big(p_B \le {\hat p}_A \le    p_B + \frac{K\sqrt{{\hat p}_A (1-{\hat p}_A)}-\varepsilon}{\sqrt{n\, \nu_n}} \big) P\big({\hat p}_B \in ( p_B - \frac{\varepsilon}{\sqrt{n\,\nu_n}}, p_B) \big).
\end{align*}
Taking logs  and limits in the above product, we have to consider two parts. For the first, we have
by Lemma \ref{lem:sqrt} below
\[
\lim_n \frac{1}{n\, \nu_n } \log P\big(p_B \le {\hat p}_A \le    p_B + \frac{K\sqrt{{\hat p}_A (1-{\hat p}_A)}-\varepsilon}{\sqrt{n\, \nu_n}} \big) = C.
\]
for some constant C; therefore,
\[
\lim_n \frac{1}{n} \log  P\big(p_B \le {\hat p}_A \le    p_B + \frac{K\sqrt{{\hat p}_A (1-{\hat p}_A)}-\varepsilon}{\sqrt{n\, \nu_n}} \big) = 0.
\]
The limit of the log of the second part divided by $n$ is 0, since
\[
P\big({\hat p}_B \in ( p_B - \frac{\varepsilon}{\sqrt{n\,\nu_n}}, p_B) \big) \ge P( -\varepsilon \le \sqrt{n(1-\nu_n)}({\hat p}_B-p_B) \le 0) \rightarrow C'>0
\]
by the CLT. \qed

\begin{lemma}\label{lem:sqrt}
Let $V_1 ,V_2 ,\ldots$ be i.i.d with $E V_1 < 0$ and moment generation function $M(t)$, and let $X_n$ be positive and uniformly bounded random variables that satisfy $X_n \stackrel{a.s.}{\rightarrow} K$ for a constant $K >0$; then,
\[
\lim_n \frac{1}{n} \log P\big( 0\le {\bar V}_n \le  \frac{X_n}{\sqrt{n}} \big) = \inf_{t>0}\{ \log M(t)  \} .
\]
\end{lemma}
\noindent {\bf Proof of Lemma \ref{lem:sqrt}.} The lemma follows by the same argument as in \citet{van der Vaart}, p. 206 (replacing $\varepsilon$ in that proof by $\frac{\tilde {K}}{\sqrt{n}}$, where $\tilde {K}$ is the bound of $X_n$); see also the proof of parts I and II of Theorem \ref{thm:test}, where a similar argument is used. \qed

\noindent{\bf Proof of Theorem \ref{nu_star}.} We will prove part I; the proof of Part II is similar. Consider the sequence of allocations $\nu'_n=\frac{\lfloor n \cdot \nu^* \rfloor}{n}$; Theorem \ref{thm:test} implies that
\begin{equation}\label{eq:lim_nu}
\lim_n \frac{1}{n} \log \left\{ 1- P\left(\frac{{\hat p}_B(n(1- \nu'_n)) -  {\hat p}_A(n \, \nu'_n)}{\sqrt{V({\hat p}_A,{\hat p}_B,n,\nu'_n)}}>K\right)\right\} = g(\nu^*).
\end{equation}

Now, let $\tilde{\nu} \neq \nu^*$, $0 \le \tilde{\nu} \le 1$,  be a limit of a certain subsequence $\{n_k\}$, i.e., $\nu^{*(1)}_{n_k} \rightarrow \tilde{\nu}$, and define $\varepsilon = (g(\tilde{\nu}) - g(\nu^*))/2$. By (\ref{eq:lim_nu}), there exists $N$ such that for $n \geq N$
\[
\frac{1}{n} \log \left\{ 1- P\left(\frac{{\hat p}_B(n(1- \nu'_n)) -  {\hat p}_A(n \, \nu'_n)}{\sqrt{V({\hat p}_A,{\hat p}_B,n,\nu'_n)}}>K\right)\right\}  < g(\nu^*) + \varepsilon.
\]
Since  $\nu^{*(1)}_{n_k} \rightarrow \tilde{\nu}$ we have by Theorem \ref{thm:test}, for large enough $k$
\[
\frac{1}{n_k} \log \left \{1- P\left(\frac{{\hat p}_B(n_k (1-\nu^{*(1)}_{n_k})) -  {\hat p}_A(n_k\, \nu^{*(1)}_{n_k})}{\sqrt{V({\hat p}_A,{\hat p}_B,n_k,\nu^{*(1)}_{n_k})}}>K\right) \right \} > g(\tilde \nu) - \varepsilon = g(\nu^*) + \varepsilon,
\]
where $g(0)=g(1)=0$. Hence, there exists $n_k > N$ such that
\[
P\left(\frac{{\hat p}_B(n_k (1-\nu^{*(1)}_{n_k})) -  {\hat p}_A(n_k\, \nu^{*(1)}_{n_k})}{\sqrt{V({\hat p}_A,{\hat p}_B,n_k,\nu^{*(1)}_{n_k})}}>K\right) <  P\left(\frac{{\hat p}_B(n_k(1- \nu'_{n_k})) -  {\hat p}_A(n_k \, \nu'_{n_k})}{\sqrt{V({\hat p}_A,{\hat p}_B,n_k,\nu'_{n_k})}}>K\right)
\]
in contradiction to the optimality of $\nu^{*(1)}_{n_k}$; therefore the limit of every converging sub-sequence is $\nu^*$. $\framebox{}$

\noindent{\bf The proofs of Theorems \ref{large_deviation} and \ref{large_deviation_gen} are omitted because they are very similar to the proofs of Theorems \ref{thm:test} and \ref{nu_star}.}

\end{document}